\documentclass[11pt]{article} 
\usepackage{amsfonts,amsmath,latexsym,amssymb,mathrsfs,amsthm,comment}

\usepackage{caption}
\usepackage{graphicx}
\usepackage{xcolor}

\let\OLDthebibliography\thebibliography
\renewcommand\thebibliography[1]{
  \OLDthebibliography{#1}
  \setlength{\parskip}{1pt}
  \setlength{\itemsep}{0pt plus 0.0ex}
}

\def\numberlikeadb{\global\def\theequation{\thesection.\arabic{equation}}}
\numberlikeadb
\newtheorem{theorem}{Theorem}[section]

\usepackage{lscape}
\usepackage{caption}
\usepackage{multirow}
\allowdisplaybreaks
\begin{document}

\title{A probabilistic proof of some integral formulas involving incomplete gamma functions
}
\author{Robert E. Gaunt\footnote{Department of Mathematics, The University of Manchester, Oxford Road, Manchester M13 9PL, UK, robert.gaunt@manchester.ac.uk}}

\date{} 
\maketitle

\vspace{-10mm}

\begin{abstract}The theory of normal variance mixture distributions is used to provide elementary derivations of closed-form expressions for the definite integrals $\int_0^\infty x^{-2\nu}\cos(bx)\gamma(\nu,\alpha x^2)\,\mathrm{d}x$ (for $\nu>1/2$, $b>0$ $\alpha>0$) and $\int_0^\infty x^{2\nu-1}\cos(bx)\Gamma(-\nu,\alpha x^2)\,\mathrm{d}x$ (for $\nu>0$, $b>0$ $\alpha>0$), where $\gamma(a,x)$ and $\Gamma(a,x)$ are the lower and upper incomplete gamma functions, respectively. The method of proof is of independent interest and could be used to derive further new definite integral formulas.
\end{abstract}


\noindent{{{\bf{AMS 2010 Subject Classification:}}} Primary  	33B20; Secondary 60E05

\section{Introduction}

The evaluation of definite integrals is a central topic in mathematics with applications throughout the mathematical sciences. Methods for evaluation of integrals and connections to other areas of mathematics can be found, for example, in the book \cite{bm04}, and many important definite integrals are collected in standard references such as \cite{gradshetyn,olver}. For $x>0$, let $\gamma(a,x)=\int_0^x t^{a-1}\mathrm{e}^{-t}\,\mathrm{d}t$ ($a>0$) and $\Gamma(a,x)=\int_x^\infty t^{a-1}\mathrm{e}^{-t}\,\mathrm{d}t$ ($a\in\mathbb{R}$) denote the lower and upper incomplete gamma functions, respectively. 
In this note, we prove the following definite integral formulas involving incomplete gamma functions.
\begin{theorem}\label{thm1}Let $\nu>1/2$, $b>0$ and $\alpha>0$. Then
\begin{align}\label{int1}\int_0^\infty x^{-2\nu}\cos(bx)\gamma(\nu,\alpha x^2)\,\mathrm{d}x=\frac{\sqrt{\pi}}{2^{2\nu}}b^{2\nu-1}\Gamma\bigg(\frac{1}{2}-\nu,\frac{b^2}{4\alpha}\bigg),
\end{align}
Now suppose $\nu>0$, $b>0$ and $\alpha>0$. Then
\begin{align}\label{int2}\int_0^\infty x^{2\nu-1}\cos(bx)\Gamma\bigg(\frac{1}{2}-\nu,\alpha x^2\bigg)\,\mathrm{d}x=\frac{2^{2\nu-1}\sqrt{\pi}}{b^{2\nu}}\gamma\bigg(\nu,\frac{b^2}{4\alpha}\bigg).
\end{align}
\end{theorem}
Formulas (\ref{int1}) and (\ref{int2}) are not known to \emph{Mathematica} nor \emph{Maple}. The integral formula (\ref{int1}) can be deduced from a suitable specialisation of formula 7.642 of \cite{gradshetyn}, although locating this formula in this manner via the standard reference \cite{gradshetyn} is not straightforward (it is not available in the section on gamma and incomplete gamma functions, and follows from a specialisation of an integral formula involving a confluent hypergeometric function). The author was unable to locate formula (\ref{int2}) in the literature. We note that the special case of formulas (\ref{int1}) and (\ref{int2}) in which $b=0$ yields the following definite integral formulas: for $\nu>1/2$ and $\alpha>0$,
\begin{align*}\int_0^\infty x^{-2\nu}\gamma(\nu,\alpha x^2)\,\mathrm{d}x=\frac{\sqrt{\pi}}{2\nu-1}\alpha^{\nu-1/2},
\end{align*}
whilst, for $\nu>0$ and $\alpha>0$,
\begin{align*}\int_0^\infty x^{2\nu-1}\Gamma\bigg(\frac{1}{2}-\nu,\alpha x^2\bigg)\,\mathrm{d}x=\frac{\sqrt{\pi}}{2 \nu}\alpha^{-\nu}.
\end{align*}
These formulas can be obtaining by letting $b\downarrow0$ in (\ref{int1}) and (\ref{int2}) and using the limiting forms $\Gamma(a,x)\sim -x^a/a$ ($a<0$, $x\downarrow0$) and $\gamma(a,x)\sim x^a/a$ ($a>0$, $x\downarrow0$) (see the series expansions 8.7.1 and 8.7.3 of \cite{olver}). Alternatively, these formulas can be derived by applying a simple change of variables to the definite integral formulas 8.14.3 and 8.14.4 of \cite{olver}.} 

In this note, we apply the theory of normal variance mixture distributions to provide an elementary proof of the integral formulas (\ref{int1}) and (\ref{int2}). The method of proof is of interest, because it offers the possibility to derive new definite integral formulas; see Section \ref{rm1} for a discussion. We remark that different probabilistic methods have been used to derive new integral formulas involving the Meijer $G$-function by \cite{gaunt}. The probabilistic method used in this present note involves generating functions, in particular moment generating functions and characteristic functions. Moment generating functions have also found application in the study the probabilistic extensions of Bernoulli polynomials, Euler polynomials and degenerate Stirling polynomials of the second kind, degenerate Bell polynomials and discrete degenerate Bell distributions \cite{chen,kim1,kim2,kim3,kim4}.

\section{Normal variance mixture distributions}

\subsection{Basic properties}

Recall that the normal distribution with mean $\mu\in\mathbb{R}$ and variance $\sigma^2\in(0,\infty)$, denoted by $N(\mu,\sigma^2)$, has probability density function (PDF) $p(x)=(2\pi\sigma^2)^{-1/2}\exp(-(x-\mu)^2/(2\sigma^2))$, $x\in\mathbb{R}$. Now let $Z\sim N(0,1)$ be a standard normal random variable with PDF $p_Z(x)=(2\pi)^{-1/2}\mathrm{e}^{-x^2/2}$, $x\in\mathbb{R}$. Suppose that $W\geq0$ is a non-negative random variable independent of $Z$. Then, the random variable
\begin{equation*}X=\sqrt{W}Z
\end{equation*} 
is said to have a normal variance mixture distribution \cite{bks82}. Normal variance mixture distributions are widely used in financial modelling, and several important distributions can be represented as normal variance mixtures, such as Student's $t$-distribution (mixing with an inverse gamma random variable) and the variance-gamma distribution (mixing with a gamma random variable); see \cite[Section 5]{fs06}. 


We now note some basic properties of normal variance mixture distributions that will be required in the sequel. Conditional on $W$, the random variable $X$ follows the $N(0,W)$ distribution.  Provided $\mathbb{P}(W=0)=0$ (that is the random variable $W$ takes the value 0 with zero probability), it thus readily follows by a simple conditioning argument that the PDF of $X$ is given by
\begin{equation}\label{pdf00}p_X(x)=\mathbb{E}\bigg[\frac{1}{\sqrt{2\pi W}}\mathrm{e}^{-x^2/(2W)}\bigg], \quad x\in\mathbb{R}.
\end{equation}
In the case that $W$ is a continuous random variable (for which the condition $\mathbb{P}(W=0)$ is automatically satisfied), the expectation in equation (\ref{pdf00}) can be expressed in terms of an integral as follows:
\begin{equation}\label{pdf}
p_X(x)=\int_{-\infty}^\infty\frac{1}{\sqrt{2\pi w}}e^{-x^2/(2w)}p_W(x)\,\mathrm{d}x.    
\end{equation}
Also, by a simple conditioning argument and using the standard fact that the characteristic function of $Y\sim N(0,w)$, $w>0$, is given by $\varphi_Y(x)=\mathbb{E}[\mathrm{e}^{\mathrm{i}tY}]=\mathrm{e}^{-t^2w/2}$, $t\in\mathbb{R}$, it follows that the characteristic function of $X$ is given by
\begin{equation}\label{cf}\varphi_X(t)=\mathbb{E}[\mathrm{e}^{\mathrm{i}tX}]=M_W(-t^2/2), \quad t\in\mathbb{R},
\end{equation}
where $M_W(u)=\mathbb{E}[\mathrm{e}^{uW}]$ is the moment generating function of $W$ (which always exists for $u\leq0$).


\subsection{Mixing with Pareto and beta distributions}

Mixing with Pareto and beta distributions was considered by \cite{r79}. Let $W_1$ follow the Pareto distribution with shape parameter $\lambda>0$ and PDF $p_{W_1}(x)=\lambda x^{-(\lambda+1)}$, $x\geq1$. Let $Z\sim N(0,1)$ be independent of $W_1$. Then the PDF of $X_1=\sqrt{W_1}Z$
is readily obtained using equation (\ref{pdf}). For $x\in\mathbb{R}$,
\begin{align}p_{X_1}(x)&=\int_1^\infty \frac{1}{\sqrt{2\pi w}}\mathrm{e}^{-x^2/(2w)}\cdot \lambda w^{-(\lambda+1)}\,\mathrm{d}w= \frac{\lambda}{\sqrt{2\pi}}\bigg(\frac{x^2}{2}\bigg)^{-\lambda-1/2}\int_0^{x^2/2}y^{\lambda-1/2}\mathrm{e}^{-y}\,\mathrm{d}y\nonumber\\
\label{pdf2}&=\frac{2^\lambda\lambda}{\sqrt{\pi}}|x|^{-2\lambda-1}\gamma\bigg(\lambda+\frac{1}{2},\frac{x^2}{2}\bigg),
\end{align}
where we used the substitution $y=x^2/(2w)$. For $u<0$, the moment generating function of $W_1$ is easily calculated in terms of the upper incomplete gamma function:
\begin{align}\label{bvf}M_{W_1}(u)=\int_1^\infty \mathrm{e}^{ux}\cdot \lambda x^{-(\lambda+1)}\,\mathrm{d}x=\lambda(-u)^\lambda\Gamma(-\lambda,-u).
\end{align}
Combining equations (\ref{cf}) and (\ref{bvf}) now yields the following formula for the characteristic function of $X_1$:
\begin{equation}\label{mgf2}\varphi_{X_1}(t)=\frac{t^{2\lambda}\lambda }{2^\lambda}\Gamma\bigg(-\lambda,\frac{t^2}{2}\bigg), \quad t\in\mathbb{R}.
\end{equation}

Now, let $W_2$ follow the beta distribution with PDF $p_{W_2}(x)=\lambda x^{\lambda-1}$, $0<x<1$, where $\lambda>0$. Suppose $Z\sim N(0,1)$ is independent of $W_2$, and let $X_2=\sqrt{W_2}Z$. By similar arguments as in the case of Pareto mixing, we obtain the following formulas for the PDF and characteristic function of the random variable $X_2$:
\begin{equation}\label{pdf3}p_{X_2}(x)=\frac{2^{-\lambda} \lambda}{\sqrt{\pi}}|x|^{2\lambda-1}\Gamma\bigg(\frac{1}{2}-\lambda,\frac{x^2}{2}\bigg), \quad x\in\mathbb{R},
\end{equation}
and
\begin{equation}\label{mgf3}\varphi_{X_2}(t)=\frac{2^{\lambda}\lambda }{t^{2\lambda}}\gamma\bigg(\lambda,\frac{t^2}{2}\bigg), \quad t\in\mathbb{R}.
\end{equation}


\section{Proof of Theorem \ref{thm1}}


\noindent{\emph{Proof of Theorem \ref{thm1}.}} 
Let us first prove the integral formula (\ref{int1}). Writing out the characteristic function $\varphi_{X_1}(t)=\int_{-\infty}^\infty \mathrm{e}^{\mathrm{i}tx}p_{X_1}(x)\,\mathrm{d}x$ as an integral, together with the formulas (\ref{pdf2}) and (\ref{mgf2}), we obtain the formula
\begin{align*}\int_{-\infty}^\infty \mathrm{e}^{\mathrm{i}tx}\cdot\frac{2^\lambda\lambda}{\sqrt{\pi}}|x|^{-2\lambda-1}\gamma\bigg(\lambda+\frac{1}{2},\frac{x^2}{2}\bigg)\,\mathrm{d}x=\frac{t^{2\lambda}\lambda }{2^\lambda}\Gamma\bigg(-\lambda,\frac{t^2}{2}\bigg).
\end{align*}
Taking the real part of the above integral and using the fact that the integrand is an even function, we obtain that
\begin{align}\label{bvfs}2\int_{0}^\infty \cos(tx)\cdot\frac{2^\lambda\lambda}{\sqrt{\pi}}x^{-2\lambda-1}\gamma\bigg(\lambda+\frac{1}{2},\frac{x^2}{2}\bigg)\,\mathrm{d}x=\frac{t^{2\lambda}\lambda }{2^\lambda}\Gamma\bigg(-\lambda,\frac{t^2}{2}\bigg).
\end{align}
We now obtain the integral formula (\ref{int1}) by making the substitution $x=\sqrt{2\alpha} y$ in the integral in (\ref{bvfs}) and then setting $\nu=\lambda+1/2$ and $b=t\sqrt{2\alpha}$.

The proof of the integral formula (\ref{int2}) is similar; we just use formulas (\ref{pdf3}) and (\ref{mgf3}) for the normal variance mixture distribution with beta mixing distribution instead of the corresponding formulas for the Pareto mixture. \hfill $\Box$

\section{Conclusion}\label{rm1}
In this note, we have used the theory of normal variance mixture distributions to provide elementary derivations of closed-form expressions for two definite integrals involving the lower and upper incomplete that are not known to \emph{Mathematica} nor \emph{Maple} and cannot be easily located in standard references on integrals.
Moreover, our method of derivation is of independent interest. Indeed, the approach of this note to obtain the integral formulas of Theorem \ref{thm1} could be used to derive integral formulas for other special functions by taking different mixture distributions beyond the Pareto and beta. The aim would be to find distributions for which the PDF and characteristic function of the normal variance mixture distribution can be calculated exactly in terms of special functions, and for which the formula for the PDF is sufficiently `complicated' so that the integral has not previously been calculated in the literature. 

We remark that the notion of a normal variance mixture distribution generalises to normal variance-mean mixture distributions \cite{bks82}. For a non-negative random variable $W$ independent of $Z\sim N(0,1)$, and constants $\mu\in\mathbb{R}$, $\theta\in\mathbb{R}$ and $\sigma>0$, we define the normal variance-mean mixture $X=\mu+\theta W+\sigma \sqrt{W}Z$. Formulas (\ref{pdf}) and (\ref{cf}) generalise in a straightforward manner for such distributions, and offer further possibility for the discovery of new integral formulas.

\section*{Acknowledgements}
I would like to thank the reviewers for their helpful comments and suggestions. The author is funded in part by EPSRC grant EP/Y008650/1.

\end{document}